\documentclass{amsart}
\numberwithin{equation}{section}

\newtheorem{theorem}{Theorem}[section]
\newtheorem{lemma}[theorem]{Lemma}

\newtheorem{definition}[theorem]{Definition}
\newtheorem{remark}[theorem]{Remark}

\begin{document}

\author{Ravshan Ashurov}
\address{Ashurov R: Institute of Mathematics of Uzbekistan,
Tashkent, Uzbekistan}
\email{ashurovr@gmail.com}

\author{Ilyoskhuja Sulaymonov}
\address{Sulaymonov I: National University of Uzbekistan,  Tashkent, Uzbekistan}
\email{ilyosxojasulaymonov@gmail.com}

\small

\title[Determining the order of ...] {Determining the order of time and spatial fractional derivatives}

\begin{abstract}

The paper considers the initial-boundary value problem for equation $D^\rho_t u(x,t)+ (-\Delta)^\sigma u(x,t)=0$, $\rho\in (0,1)$, $\sigma>0$, in an N-dimensional domain $\Omega$ with a homogeneous Dirichlet condition. The fractional derivative is taken in the sense of Caputo. 
	The main goal of the work is to solve the inverse problem of simultaneously determining two parameters: the order of the fractional derivative $\rho$ and the degree of the Laplace operator $\sigma$. A new formulation and solution method for this inverse problem are proposed. It is proved that in the new formulation the solution to the inverse problem exists and is unique for an arbitrary initial function from the class $L_2(\Omega)$. Note that in previously known works, only the uniqueness of the solution to the inverse problem was proved and the initial function was required to be sufficiently smooth and non-negative.

\end{abstract}

\maketitle

\section{Introduction}\label{s:1}

Let $\Omega \subset \mathbb{R}^N$ be an arbitrary bounded domain, with sufficiently smooth boundary $\partial\Omega$, and $\rho\in (0,1)$, $\sigma>0$. Consider the initial-boundary value problem for a space-time-fractional diffusion equation
\begin{equation}\label{laplas}
	\left\{
	\begin{aligned}
		&D^\rho_t u(x,t)+(-\Delta)^\sigma u(x,t)=0, \quad x\in \Omega, \quad 0<t\leq T, \\
		&u(x, t)=0, \quad x\in \partial \Omega, \quad 0<t\leq T, \\
		&u(x,0)=\varphi(x), \quad x\in \Omega, \\
	\end{aligned}
	\right.
\end{equation}
where $\Delta$ is the Laplace operator, the fractional power of this operator is determined using the von Neumann spectral theorem, $\varphi(x)\in L_2(\Omega)$, $D^\rho$ is the fractional Caputo derivative of order $0<\rho< 1$, which for any locally integrable function $h(t)$ is defined as (see e.g. \cite{Kilbas}, p. 91)
\begin{equation}\label{HT}
	D_t^\rho h(t)=\frac{1}{\Gamma
		(1-\rho)}\frac{d}{dt}\int\limits_0^t\frac{h(\xi)-h(0)}{(t-\xi)^{\rho}} d\xi,
	\quad t>0,
\end{equation}
provided the right-hand side exists.  Here $\Gamma(\rho)$ is
Euler's gamma function. We call problem (\ref{laplas}) \emph{the forward problem}.

Even in the middle of the last century in his visionary papers \cite{Riesz1} and \cite{Riesz2} Marcel Riesz introduced the fractional
powers of the Laplacean and string vibration operator in Euclidean and Lorentzian spaces and developed the calculus of these non-local operators and studied the Dirichlet and Cauchy problems for
respectively $(-\Delta)^\sigma$ and $(\partial_{tt} - \Delta)^\sigma$. These pseudo-differential operators play an important role in many branches of the applied sciences ranging from fluid dynamics to elasticity
and to quantum mechanics (see, e.g., \cite{ManyDefinitions}). 

The spatial fractional diffusion equations have recently attracted considerable attention due to their applications in a soft matter (or soft condensed matter), elasticity, turbulence, anomalous diffusion, finance, biology, noise reduction, etc (see, e.g. \cite{Chen, Dipierro, Oleg, woj, Meer1, Hilfer, Meer2}). The space-time fractional diffusion equation (\ref{laplas}) is also related to the limitations of the scale of random walks in continuous time (see \cite{Meer1, Hilfer, Meer2}).

Based on the above applications, it is undoubtedly topical to find a solution to the formulated problem (\ref{laplas}). Works \cite{Chen1} and \cite{Yamamoto} (in this paper, instead of the Laplace opretor, the author considered a second-order elliptic operator with smooth coefficients) are devoted to the study of precisely this question. The authors of these papers wrote out the strong solution of the forward problem (\ref{laplas}) using Fourier series.

In work \cite{AZ}, a problem similar to (\ref{laplas}) is considered in the whole space $\mathbb{R}^N$. In this paper, the Laplace operator is replaced by an elliptic operator of arbitrary order with constant coefficients. Using the Fourier transform, the authors found a classical solution to the corresponding initial-boundary value problem.

The purpose of this work is to solve the inverse problem of determining both the order of the fractional derivative $\rho$ and the degree $\sigma$ of the Laplace operator in the intial-boundary value problem (\ref{laplas}).

When dealing with fractional order equations, it is often challenging to determine the order of the fractional derivative directly because there is no measurement instrument available. In such cases, it becomes necessary to approach the inverse problem, which involves inferring this parameter by utilizing indirectly observable information about the solutions. These inverse problems hold not only theoretical interest but also practical significance in solving initial-boundary value problems and exploring solution properties. An extensive selection of articles published up to 2019 can be found at the link provided in the review article Z. Li, Y. Liu, M. Yamamoto \cite{Li}. We will not dwell separately on each work cited here, but will only make the following general remark. In all works, the unknown parameter is found using a given value of the solution to the forward problem over a long period of time, i.e. an additional condition for solving the inverse problem looks like this:
\begin{equation}\label{AdYamamoto}
	u(x_0, t)=h(t), \,\, x_0\in \Omega,\,\, 0<t<T.    
\end{equation}
The authors managed to prove only the uniqueness of the solution to the inverse problem (with the exception of work J. Janno \cite{Jan}, where the existence was also proven). The solution methods are based on the asymptotic behavior of the Mittag-Leffler functions at infinity. Namely, by the method of separation of variables, two solutions $u_1(x,t)$ and $u_2(x,t)$ to the forward problem  are constructed corresponding to two different values of the order of fractional derivative $\rho_1$ and $\rho_2$. The equality of two solutions for $x=x_0\in \Omega$ and $t\in (0, T)$, due to the asymptotics of the Mittag-Leffler functions, gives the equality of two parameters $\rho_1$ and $\rho_2$. 

The same method was used in a recent paper X. Jing and M. Yamamoto \cite{JY}, which considers the fractional diffusion-wave equation in a one-dimensional spatial domain. The authors studied the inverse problem by simultaneously determining three parameters of the problem: the order of the fractional derivative, the potential (i.e., the coefficient of $u(x,t)$) and the initial functions. It has been proven that if the solution to the direct problem is given at two points $x_1, x_2$ and $t\in (0,T)$, then the solution to the inverse problem is unique.

Of course, the uniqueness of the solution to the inverse problem is very important from the application point of view, since when finding parameter $\rho$ using a numerical method, we must be sure that we have found exactly the solution we need. The uniqueness of the solution guarantees this. But on the other hand, if a solution to the inverse problem does not exist, this means that we have built an incorrect model equation for the process under study. 

In this regard, a natural question arises: is there a solution to the inverse problem for an arbitrarily given function $h(t)$ in the additional condition (\ref{AdYamamoto})? The answer is no. Since point $x_0$ is located inside the domain where the equation is considered, then function $h(t)$ must satisfy the subdiffusion equation and since any solution to such an equation is infinitely differentiable with respect to $t$, then function $h(t)$ must be at least infinitely differentiable with respect to $t$. And this is certainly not enough; the complexity of the proof of
the existence can be seen from the statement of corresponding Theorem 7.2 of J. Janno \cite{Jan}, which is formulated on more than one journal page.

Another motivation in favor of considering another additional condition than (\ref{AdYamamoto}) for the inverse problem is expressed in the above review work Z. Li, Y. Liu, M. Yamamoto \cite{Li}. In the section  Open Problems the authors noted: "The studies on inverse problems of the recovery of the fractional orders are far from satisfactory since all the publications ... studied this inverse problem by measuring over $t \in (0, \infty)$. It would be interesting to investigate inverse problem by the value of the solution at a fixed time as the observation data." Also note that it would be mathematically correct to specify not a function, as in condition (\ref{AdYamamoto}), to search for one unknown number $\rho$, but to specify one number.

Based on these considerations, other additional conditions were proposed in works \cite{AU, AA1, AA2, AF12, AZ10, AF13, AU1} to solve this inverse problem. Thus, in the work \cite{AU} it was proven that the solution to the inverse problem with the additional condition $\int_\Omega u(x, t_0) v_1(x) dx=d_0$ exists and is unique, where $v_1(x)$ is the first eigenfunction of the Laplace operator with the Dirichlet condition, $t_0$ is a sufficiently large number and $d_0$ is a given number. A similar result is established in works \cite{AA1} and \cite{AA2} with an additional condition $\int_\Omega |u(x, t_0)|^2 dx = d_0$. In works \cite{AZ10, AF13, AU1} inverse problems are considered, respectively, for equations of mixed type, for a fractional wave-equation and for systems of fractional pseudo-differential equations. The authors of the work \cite{AF12} investigated the inverse problem of simultaneously determining the order of the fractional derivative and the right-hand side of the subdiffusion equation. In all the listed works \cite{AU, AA1, AA2, AF12, AZ10, AF13, AU1}, the solution to the inverse problem is based on the monotonicity of the Mittag-Leffler functions $E_\rho(t)$ with respect to the parameter $\rho$. Note that such monotonicity was first discovered in these listed works. Since our attention is focused on studying the inverse problem of simultaneous determination of two parameters, we will not go into detail in these works, but will focus specifically on those works where such two-parameter problems have been studied.

To the best of our knowledge, only the next three papers \cite{Yamamoto} \cite{AZ} and \cite{TU} studied inverse problems for the simultaneous determination of two parameters. Moreover, in \cite{Yamamoto} and \cite{TU}, only uniqueness theorems were proved. And in \cite{AZ}, the authors succeeded, by setting the over-determination condition in a different form, to prove both the existence and uniqueness of the desired parameter. Let us consider these works in more detail.

Tatar and Ulusoy \cite{TU} considered the following initial-boundary value
problem 
\[
\partial_t^\rho u(x,t) +  (-\Delta)^\sigma u(x, t) =0, \quad 0<t<T, \ x \in
(-1,1),
\]
\[
u(-1,t)=u(1,t)=0,\quad 0<t<T,
\]
\[
u(x,0)=\varphi(x), \quad x \in
(-1,1),
\]
where $(-\Delta)^\sigma$ is the one-dimensional fractional Laplace
operator, $\rho \in (0,1)$ and $\sigma \in (1/4,1)$. The fact that the space dimension here is equal to $1$ seems to be significant, since when solving this problem by the Fourier method, the authors use the asymptotics of the eigenvalues of the one-dimensional Laplacian. The main result of the work is the proof of the uniqueness of the solution of the inverse problem for the determination of two parameters $\rho$ and $\sigma$ with an additional condition $u(0,t)=h(t)$, $0<t<T$. When proving the uniqueness theorem, the authors had to set the following rather stringent conditions on the initial function:
$$
\varphi \in C^\infty (-1,1), \,\, \varphi_k> 0, \text{for
	all}\quad k\geq 1,
$$
where $\varphi_k$ are the Fourier coefficients of $\varphi$.

Later, in 2020, a similar two parameter inverse problem in an $N$-dimensional bounded domain $\Omega$ with smooth boundary $\partial \Omega$ was studied by M. Yamamoto \cite{Yamamoto} with the same additional condition and $\sigma \in (0,1)$. He proved a uniqueness theorem with a less stringent condition on the initial function: $\varphi$ is equal to zero on $\partial \Omega$, $\varphi \in L_2^\tau(\Omega)$, $\tau>N/2$ and $\varphi(x)\geq 0 $ for all $x\in\Omega$.
Here $L_2^\tau(\Omega)$ denotes the Sobolev space. Note that to prove the uniqueness of the solution to the inverse problem, the author used the above-mentioned method based on the asymptotic behavior of the Mittag-Leffler functions.

In work \cite{AZ}, a two parameter inverse problem is studied in the case when $\Omega= \mathbb{R}^N$ and instead of the Laplace operator, an arbitrary elliptic operator with constant coefficients is considered. The authors of this work managed to find such additional conditions that guarantee not only the uniqueness but also the existence of a solution to the inverse problem. It should also be noted that in \cite{AZ} no additional conditions are imposed on the initial function $\varphi$, except for smoothness $\varphi \in L_2^\tau(\mathbb{R}^N)$, $\tau>N/2$. It should be emphasized that in the proof of this result it is important that the elliptic operator in question has a continuous spectrum, namely the fact that the spectrum passes through the point $\lambda=1$. Analysis of the applied method shows that it can be applied for the case when elliptic part of the subdiffusion equation has the discreet spectrum only in the case when $\lambda=1$ is the eigenvalue of this operator.

In the present paper, we study the initial-boundary value problem (\ref{laplas}), in which the parameters $\rho$ and $\sigma$ are assumed to be unknown. Such an inverse problem will be called \emph{the two parameter inverse problem}. Note that the equation in (\ref{laplas}) for any positive $\sigma>0$ is hypoelliptic and therefore we allow the parameter $\sigma$ to take any positive values. It follows from the main result of this work that if the value of the function (note that the solution to the problem (\ref{laplas}) depends on both $\rho$ and $\sigma$)
\begin{equation}\label{AdConGen}
	\int_\Omega u(x,t)v_1(x)dx,
\end{equation}
is set at two different points $t_0$ and $t_1$, then the desired parameters $\rho$ and $\sigma$ can be uniquely determined.
Here the symbol $v_1(x)$ denotes the first eigenfunction of the Laplace operator with the Dirichlet condition.

If in problem (\ref{laplas}) the parameter $\sigma$ is known and only the solution $u(x,t)$ of the initial-boundary value problem and the order of the fractional derivative $\rho$ are unknown, then we will call such an inverse problem \emph{ the first inverse problem}. It is also proved in the paper that for the unique restoration of the parameter $\rho$ it is sufficient to know the value of the function (\ref{AdConGen}) only at one point $t_0$.

This article is divided into six sections.
In the next section, we give exact statements of the problems under study and formulate the results of the work. In Section \ref{s:3}, we study some properties of Mittag-Leffler function. In Section \ref{s:4}, we study the forward problem (\ref{laplas}). Section \ref{s:5} is devoted to the study of the first and two parameter inverse problems. In the last section conclusion is presented.

\section{Problem statement and results}\label{c:2}

Let us start with how we define the fractional power of the Laplace operator. Over the years, researchers have explored the fractional Laplacian from many perspectives, such as a probabilistic approach, potential theory, or partial differential equations. Different views lead to different definitions of the fractional Laplacian. If we consider the Laplacian in $\mathbb{R}^N$, then all known definitions of fractional powers turn out to be equivalent under certain assumptions \cite{10Definitions}. However, in a bounded domain, which is crucial for physical modeling, there is no single definition of the fractional Laplacian \cite{ManyDefinitions}. There are basically three definitions in the literature: the restricted fractional Laplacian, the spectral fractional Laplacian (based on von Neumann's theorem), and the regional Laplacian (see, for example, \cite{ManyDefinitions} and \cite{10Definitions}).

In this paper, the fractional Laplacian is defined using the von Neumann spectral theorem. This definition uses only eigenfunctions and eigenvalues of the Laplace operator with the Direchlet condition. Therefore, instead of the Laplace operator, one can consider an arbitrary positive self-adjoint operator.

Let $H$ be a separable Hilbert space with the scalar product $(\cdot,\cdot)$ and the norm $||\cdot||$. Consider the operator $A$ satisfying the following conditions
\begin{enumerate}
	\item $(Ah,h)\geq C(h,h), h\in D(A)$, $C>0$,
	\item $A=A^*$.
\end{enumerate}
Here $D(A)$ is the definition domain of operator $A$ (given below) and $A^*$ is the adjoint of operator $A$.

Let $A$ be an operator with complete orthonormal eigenfunctions $v_k$ in $H$ and countable and ordered eigenvalues $\lambda_k: 0<\lambda_1\leq \lambda_2\leq \dots \rightarrow+\infty$ . Additionally, we assume that the sequence ${\lambda_k}$ does not have any finite limit points.

Let $\sigma$ be a finite real number. We introduce the power of operator $A$, acting in $H$
according to the rule:
\begin{equation}\label{eq1:1}
	A^\sigma g =\sum\limits_{k=1}^\infty \lambda_k^\sigma g_kv_k,
\end{equation}
where $g_k=(g,v_k)$ are the Fourier coefficients of the vector $g\in H$. The domain of this operator is given by the following expression:
$$
D(A^\sigma)=\{g\in H:  \sum\limits_{k=1}^\infty \lambda_k^{2\sigma}
|g_k|^2 < \infty\}.
$$

For a vector-valued function (or simply a function) $h: \mathbb{R}_+\rightarrow H$, the Caputo fractional derivative of order $0<\rho<1$ is defined as formula (\ref{HT}) (see \cite{Lizama}).

Consider the following problem:
\begin{equation}\label{eq1:3}
	\left\{
	\begin{aligned}
		&D^\rho_t u(t)+A^\sigma u(t)=0, \quad  0<t\leq T, \\
		&u(0)=\varphi, 
	\end{aligned}
	\right.
\end{equation}
where function $\varphi \in H$ and $\rho\in (0,1), \sigma>0$. This problem is also called \emph{the forward problem}.

\begin{definition}\label{def1}
	A function \,  $u(t)\in C([0,T];H)$  \, with the properties
	\, $D_t^\rho u(t)$,\, $Au(t)\in C((0,T];H)$ and
	satisfying conditions (\ref{eq1:3})  is called \textbf{the
		solution} of the problem (\ref{eq1:3}).
\end{definition}

Now we can formulate the existence and uniqueness theorem for the forward problem.

\begin{theorem}\label{t1}
	Let $\varphi\in H$. Then the forward problem has a unique solution and this solution has the form
	\begin{equation}\label{eq1:4}
		u(t)=\sum_{k=1}^\infty \varphi_k E_{\rho}(-\lambda_k^\sigma t^\rho)v_k.
	\end{equation}
\end{theorem}
Recall that the Mittag-Leffler function $E_\rho(t)$ is defined as follows:
$$
E_\rho(t)=\sum_{k=0}^\infty \frac{t^k}{\Gamma(\rho k+1)}.
$$

\begin{remark}Note that Theorem \ref{t1} in particular implies the existence and uniqueness of a strong solution to problem (\ref{laplas}) for any initial function $\varphi\in L_2(\Omega)$.
\end{remark}

In this article we also study two inverse problems. First, assuming that the degree $\sigma$ of the operator $A$ is known and the parameter $\rho$ is unknown, we study the inverse problem of determining this parameter. As an additional condition, we take the value (\ref{AdConGen}) at $t_0$:
\begin{equation}\label{q1}
	|(u(t_0),v_1)|=d_0.
\end{equation}

As noted above, problem (\ref{eq1:3}) with the additional condition  (\ref{q1}) will be called the first inverse problem.

\begin{definition}\label{inpr1}
	A pair $\{u(t),\rho\}$ of the function $u(t)$ satisfying all the conditions of Definition \ref{def1}  and the parameter $\rho\in(0,1)$ is called the solution of the first inverse problem (\ref{eq1:3}), (\ref{q1}).
\end{definition}

To solve the first inverse problem fix the number $\rho_0\in(0,1)$ and consider the problem for $\rho\in[\rho_0,1)$.

\begin{theorem}\label{t2}
	Let $\varphi\in H$ and $\sigma>0$ be given, and $\lambda_1$ be the first eigenvalue of operator $A$. Then there is a number $T_0=T_0(\lambda_1,\rho_0)$ such that for all $t_0\geq T_0$ there exists a unique solution $\{u(t),\rho\}$ of the first inverse problem if and only if
	\begin{equation}\label{ty1}
		e^{-\lambda_1^\sigma t_0} \leq \frac{d_0}{|\varphi_1|}< E_{\rho_0}(- t_0^{\rho_0}).
	\end{equation}
\end{theorem}

Now we consider both parameters $\rho$ and $\sigma$ unknown and study the inverse problem of determining both parameters simultaneously. As the next additional condition we take the value
of (\ref{AdConGen}) at $t_1$:

\begin{equation}\label{q2}
	|(u(t_1),v_1)|=d_1. 
\end{equation}

Problem (\ref{eq1:3}) with additional conditions (\ref{q1}) and (\ref{q2}) is called \emph{the two parameter inverse problem}.

Let us pay attention to the following fact: it is easy to see that
$$
|(u(t_1),v_1)|=|\varphi_1|E_\rho(-\lambda_1^{\sigma} t_1^{\rho}).
$$
Obviously, it follows that, if $\lambda_1=1$, then in equation (\ref{q2}) to determine $\sigma$, this unknown disappears. Consequently, when studying the two parameter inverse problem in (\ref{q1}) and (\ref{q2}), instead of $v_1$ it is necessary to take that eigenfunction whose corresponding eigenvalue is not equal to 1. In this case, further reasoning does not change. Therefore, in what follows, without loss of generality, we will assume that $\lambda_1\neq 1$.

\begin{definition}\label{inpr1}
	A triple $\{u(t),\rho,\sigma\}$ of the function $u(t)$ satisfying all the conditions of Definition \ref{def1} and the parameters $\rho\in(0,1)$ and $\sigma>0$ is called the solution of the inverse problem (\ref{eq1:1}),(\ref{q1}) and (\ref{q2}).
\end{definition}

To solve the two parameter inverse problem fix the numbers $\rho_0,\rho_1\in(0,1)$ ($\rho_0<\rho_1$) and consider the problem for $\rho\in[\rho_0,\rho_1]$ and  $\sigma>0$.

\begin{theorem}\label{t3}
	Let $\lambda_1$ be the first eigenvalue of operator $A$. Then there are numbers $T_1=T_1(\lambda_1,\rho_1)>0$ and $c_0=c_0(\rho_1)>1$ such that for all $t_0>T_1$ and $t_1>T_1$ satisfying condition $t_0>c_0 t_1$, there is a unique solution $ \{u(t),\rho,\sigma\}$ of the inverse problem (\ref{eq1:3}), (\ref{q1}), (\ref{q2}), if and only if $d_0$ satisfies inequality (\ref{ty1}) and $d_1$ satisfies the following inequality 
	\begin{equation}\label{ty2}
		E_{\rho_1}(-\lambda_1^\sigma t_1^{\rho_1}) \leq \frac{d_1}{|\varphi_1|}< E_{\rho_0}(-t_1^{\rho_0}).
	\end{equation}
\end{theorem}

\section{Properties of Mittag-Leffler function}\label{s:3}

For the Mittag-Leffler function with a negative argument we have an estimate (see, for example, \cite{Gor}, p. 29)
\begin{equation}\label{eq1:2}
	|E_{\rho}(-t)|\leq \frac{C}{1+ t}, \quad t>0.
\end{equation}

Consider the function $F(\xi,\rho)$:
\begin{equation}\label{yorf1}
	F(\xi,\rho)=\frac{1}{2\pi i \lambda^\sigma  \rho t^\rho}\frac{e^{\xi^{1/\rho}}\xi}{\xi+\lambda^\sigma t^\rho},
\end{equation}
and study the following two integrals:
\begin{equation}\label{yf1}
	f_{\pm}(\rho)=e^{\pm i\beta}\int\limits_1^\infty F(se^{\pm i\beta},\rho)ds,   
\end{equation}

\begin{equation}\label{yf2}
	g(\rho)=i\beta \int_{-1}^1F(e^{i\beta s},\rho)e^{i\beta s}ds,    
\end{equation}
where $\rho\in(0,1)$ and $\beta=\frac{3\pi}{4}\rho$.

\begin{lemma}\label{lem2}
	There is a constant $C>0$, such that
	$$
	|f'_{\pm}(\rho)|\leq \frac{C}{(\lambda^\sigma t^\rho)^2}\left[\frac{1}{\rho}+\ln t\right], \quad t>1, \,\, \lambda>0.
	$$
\end{lemma}
{\em Proof.} Consider first $f_+(\rho)$. Since $\beta=\frac{3\pi}{4}\rho$ and $\xi=se^{i\beta}$, then $e^{\xi^{1/\rho}}=e ^{\frac{1}{\sqrt{2}}(i-1)s^{\frac{1}{\rho}}}$. Calculate the derivative $f'_+(\rho)$ and get:
$$
f'_+(\rho)=\frac{1}{2\pi i \lambda^\sigma  \rho t^\rho}\times
$$
$$
\times \int\limits_1^\infty \frac{e^{\frac{1}{\sqrt{2}}(i-1)s^{\frac{1}{\rho}}}se^{2ia\rho}\left[-\frac{i-1}{\sqrt{2}\rho^2}s^{1/\rho}\ln s+2ia-\frac{1}{\rho}-\ln t -\frac{iase^{ia\rho}+\lambda^\sigma t^\rho \ln t}{se^{ia\rho}+\lambda^\sigma t^\rho} \right]}{se^{ia\rho}+\lambda^\sigma t^\rho}ds,
$$
where $a=\frac{3\pi}{4}$. Since $|se^{ia\rho}+\lambda^\sigma t^\rho|\geq \lambda^\sigma t^\rho$ then:
$$ 
|f' _+(\rho)|\leq \frac{C}{\rho (\lambda^\sigma t^\rho)^2}\int\limits_1^\infty e^{-\frac{1}{\sqrt {2}}s^{1/\rho}}s\left[\frac{1}{\rho^2}s^{1/\rho}\ln s+\ln t \right]ds. 
$$
If we define $r=s^{\frac{1}{\rho}}$, then $s=r^\rho, ds=\rho r^{\rho-1}dr$ and since $\ln s^{\frac{1}{\rho}}<s^{\frac{1}{\rho}}$, we have:
$$
|f'_+(\rho)|\leq \frac{C}{\rho (\lambda^\sigma t^\rho)^2}\int\limits_1^\infty e^{-\frac{1}{\sqrt{2}}s^{1/\rho}}s\left[\frac{1}{\rho^2}s^{1/\rho}\ln s+\ln t \right]ds\leq 
$$
$$
\leq \frac{C}{ (\lambda^\sigma t^\rho)^2}\int\limits_1^\infty e^{-\frac{1}{\sqrt{2}}r}r^{2\rho-1}\left[\frac{1}{\rho}r^2+\ln t \right]dr=
$$
$$
= \frac{C}{ (\lambda^\sigma t^\rho)^2}\left(\int\limits_1^\infty e^{-\frac{1}{\sqrt{2}}r}r^{2\rho+1}\frac{1}{\rho}dr+\int\limits_1^\infty e^{-\frac{1}{\sqrt{2}}r}r^{2\rho-1}\ln t dr\right)\leq 
$$
$$
\leq \frac{C}{ (\lambda^\sigma t^\rho)^2}\left(\frac{C_2}{\rho}+C_1\ln t\right)\leq \frac{C}{ (\lambda^\sigma t^\rho)^2}\left(\frac{1}{\rho}+\ln t\right).
$$

Function $f'_-(\rho)$ is estimated in the same way.

Lemma \ref{lem2} is proved.

\begin{lemma}\label{lem3}
	There is a constant $C>0$, such that
	$$
	|g'(\rho)|\leq C\frac{\ln t}{(\lambda^\sigma t^\rho)^2}, \quad t>1, \,\,\lambda >0.
	$$
\end{lemma}
{\em Proof.} We have the following equality for the derivative $g'(\rho)$:
$$
g'(\rho)=\frac{a}{2\pi \lambda^\sigma t^\rho}\int\limits_{-1}^1\frac{e^{e^{ias}}e^{2ia\rho s}\left[ 2ias-\ln t - \frac{iase^{ia\rho s}+\lambda^\sigma t^\rho \ln t}{e^{ia\rho s}+\lambda^\sigma t^\rho}\right]}{e^{ia\rho s}+\lambda^\sigma t^\rho}ds,
$$
where $a=\frac{3\pi}{4}$. Since $|e^{ia\rho s}+\lambda^\sigma t^\rho|\geq \lambda^\sigma t^\rho$, then:
$$
|g'(\rho)|\leq C\frac{\ln t}{(\lambda^\sigma t^\rho)^2}.
$$

Lemma \ref{lem3} is proved.

Let us represent by  $\delta(1;\beta)$ a contour oriented by non-decreasing $\arg \xi$ consisting of the following parts: the ray $\arg \xi=-\beta$ with $|\xi|\geq 1$, the arc $-\beta\leq \arg\xi\leq\beta$, $|\xi|=1$, and the ray $\arg\xi=\beta$, $|\xi|\geq 1$. The counter $\delta(1;\beta)$ is called Hankel path.

Let $\beta=\frac{3\pi}{4}\rho$, $\rho\in[\rho_0,1)$ and $\sigma>0$. Then, we write function $E_{\rho}(-\lambda^\sigma t^\rho)$ in the following form ( see \cite{Gor}, p. 27):
\begin{equation}\label{eq3:1}
	E_{\rho}(-\lambda^\sigma t^\rho)=p(\rho,\sigma,t)+q(\rho,\sigma,t),
\end{equation}
where
$$
p(\rho,\sigma,t)=\frac{1}{\Gamma(1-\rho)\lambda^\sigma t^\rho },
$$
$$
q(\rho,\sigma,t)=-\frac{1}{2\pi i \lambda^\sigma  \rho t^\rho}\int_{\delta(1;\beta)}\frac{e^{\xi^{1/\rho}}\xi}{\xi+\lambda^\sigma t^\rho}d\xi.
$$

\begin{lemma}\label{lem1}
	Let $\sigma>0$ and $0<\rho_0<1$. Then there exists a number $T_0=T_0(\lambda,\rho_0)$ such that for all $t\geq T_0$ and $\lambda\geq \lambda_1$, function $E_{\rho}(-\lambda^\sigma t^\rho)$ is positive and monotonically decreasing for all $\rho\in [\rho_0,1)$ and the following estimates hold:
	\begin{equation}\label{p}
		\frac{\partial}{\partial\rho}p(\rho,\sigma,t)\leq -\frac{1}{\lambda^\sigma t^\rho},\quad t>1,\quad \lambda>0,
	\end{equation}
	\begin{equation}\label{q}
		\left|\frac{\partial}{\partial\rho}q(\rho,\sigma,t)\right|\leq C\frac{1/\rho+\ln t}{(\lambda^\sigma t^\rho)^2},\quad t>1,\quad \lambda>0.
	\end{equation}
\end{lemma}

{\em Proof.} Let us calculate the derivative of the function $p(\rho,\sigma,t)$ (see (\ref{eq3:1})) and apply the equality $\Gamma'(\rho)=\Gamma(\rho) \Psi(\rho)$, where $\Psi( \rho)$ is the logarithmic derivative of the function $\Gamma(\rho)$: $\Psi( \rho)=\frac{d}{d\rho }\ln\Gamma(\rho)$ (see \cite{Batemen}, page 15). Then we get:
$$
\frac{\partial}{\partial\rho}p(\rho,\sigma,t)=-\frac{\ln t-\Psi(1-\rho)}{\Gamma(1-\rho)\lambda^\sigma t^\rho }.
$$ 
Apply the following equalities for $\Gamma(\rho)$ and $\Psi(\rho)$:
\begin{equation}\label{eq3:2}
	\frac{1}{\Gamma(1-\rho)}=\frac{1-\rho}{\Gamma(2-\rho)}, 
\end{equation}
\begin{equation}\label{eq3:3}
	\Psi(1-\rho)=\Psi(2-\rho)-\frac{1}{1-\rho},
\end{equation}
to get
$$
\frac{\partial}{\partial\rho}p(\rho,\sigma,t)=-\frac{(1-\rho)[\ln t-\Psi(2 -\rho)]+1}{\Gamma(2-\rho)\lambda^\sigma t^\rho }.
$$ 

Since $-\gamma<\Psi(2-\rho)<1-\gamma$ (see \cite{Batemen}, p. 15) and choosing $t\geq 2$, we have the following estimate for $\frac{\partial}{\partial\rho}p(\rho,\sigma,t)$:

\begin{equation}\label{es1}
	-\frac{\partial}{\partial\rho}p(\rho,\sigma,t)\geq \frac{(1-\rho)[\ln t-(1-\gamma)]+1}{\Gamma(2-\rho)\lambda^\sigma t^\rho }\geq \frac{1}{\lambda^\sigma t^\rho}.
\end{equation}
where $\gamma\approx 0.57722$ is the Euler-Mascherano constant.

Let us calculate the derivative of the function $q(\rho,\sigma,t)$. We rewrite it in the following form:
$$
q(\rho,\sigma,t)=q_{2+}(\rho,\sigma,t)+q_{2-}(\rho,\sigma,t)+q_{21}(\rho,\sigma,t),
$$
where the functions $q_{2\pm}(\rho,\sigma,t)$ and $q_{21}(\rho,\sigma,t)$ are equal to the value of the functions (\ref{yf1}) and (\ref{yf2}), respectively.

Using Lemma \ref{lem2} and Lemma \ref{lem3}, we obtain:
\begin{equation}\label{es2}
	\left|\frac{\partial}{\partial\rho}q(\rho,\sigma,t)\right|\leq C\frac{1/\rho+\ln t}{(\lambda^\sigma t^\rho)^2}.
\end{equation}

By (\ref{es1}) and (\ref{es2}) we get:
\begin{equation}\label{es3}
	\frac{d}{d\rho}(E_{\rho}(-\lambda^\sigma t^\rho))<-\frac{1}{\lambda^\sigma t^\rho}+C\frac{1/\rho+\ln t}{\lambda^{2\sigma} t^{2\rho}}.
\end{equation}
Hence, if the following condition is met
$$
t^\rho > C\frac{1/\rho+\ln t}{\lambda^\sigma},
$$
or
$$
t^{\rho_0} > C\frac{1/\rho_0+\ln t}{\lambda^\sigma}.
$$
then $\frac{d}{d\rho}E_{\rho}(-\lambda^\sigma t^\rho)<0$ for all $\rho\in[\rho_0,1)$.
Therefore there exists a number $T_0=T_0(\lambda,\rho_0)$ such that for all $t\geq T_0$ one has:
$$
\frac{\partial}{\partial\rho}E_{\rho}(-\lambda^\sigma t^\rho)<0, \quad   \rho\in[\rho_0,1).
$$
Thus Lemma \ref{lem1} is proved. Note that estimates (\ref{p}) and (\ref{q}) were proven above (see (\ref{es1}) and (\ref{es2})).

\begin{lemma}\label{lem4}
	There is a constant $C>0$, such that
	\begin{equation}\label{es4}
		\left|\frac{\partial}{\partial\sigma}q(\rho,\sigma,t)\right|<\frac{C|\ln\lambda|}{\lambda^{2\sigma}t^{\rho}}, \quad t>1, \,\,\lambda>0\,\, (\lambda\neq 1).
	\end{equation}
\end{lemma}

{\em Proof.} Apply formula (\ref{yorf1}) to calculate derivative $\frac{\partial}{\partial\sigma}q(\rho,\sigma,t )$. Then
$$
\frac{\partial}{\partial\sigma}q(\rho,\sigma,t) = \frac{(1+t^\rho)\ln\lambda}{2\pi i \rho \lambda^\sigma t^{\rho}}\int\limits_{\delta(1;\beta)}\frac{e^{\xi^{1/\rho}}\xi}{\xi+\lambda^\sigma t^\rho} d\xi.
$$

Since the estimate of the derivative $\frac{\partial}{\partial\sigma}q(\rho,\sigma,t )$ is the similarly as the estimate of the derivative of the function $f_\pm(\rho)$ in Lemma \ref{lem2}, then we have:
$$
\left|\frac{\partial}{\partial\sigma}q(\rho,\sigma,t)\right| \leq \frac{(1+t^\rho)|\ln\lambda|}{\pi (\lambda^\sigma t^{\rho})^2}\left(C_0+\frac{4}{3}\pi\right)\leq \frac{C|\ln\lambda|}{\lambda^{2\sigma}t^\rho}.
$$

Lemma \ref{lem4} is proved.

\section{Forward problem}\label{s:4}

In this section, we study the problem of finding the solution of the forward problem $u(t)$. 

Assume that a solution of the forward problem exists. Since system $\{v_k\}$ is
complete in $H$, then this solution has the form:
\begin{equation}\label{eq2:1}
	u(t)=\sum\limits_{k=1}^\infty T_k(t)v_k,
\end{equation}
where $T_k(t)=(u(t),v_k)$ are the Fourier coefficients of the function $u(t)$ and are unknown.

Multiply the orthonormal eigenfunctions $v_k$ to the equation in the problem (\ref{eq1:3}), to get
\begin{equation}\label{eq2:2}
	(D_t^\rho u(t),v_k)+(A^\sigma u(t),v_k)=0.
\end{equation}

We have $(D_t^\rho u(t),v_k)=D_t^\rho T_k(t)$ and since operator $A$  is self-adjoint, then  $(A^\sigma u(t),v_k)=\lambda_k^\sigma T_k(t)$. Therefore, to determine $T_k(t)$ we obtain the following Cauchy problem
$$
D_t^\rho T_k(t)+\lambda_k^\sigma T_k(t)=0,\quad \quad T_k(0)=\varphi_k.
$$
This problem has a unique solution (see, for example, \cite{Gor}, p. 174):
\begin{equation}\label{eq2:3}
	T_k(t)=\varphi_k E_{\rho} (-\lambda_k^\sigma t^\rho).
\end{equation}

From this, in particular, it follows that if a solution to the forward problem
exists, then it is unique. Indeed, for this it is sufficient to prove that the solution $u(t)$ to the forward problem  with the homogeneous condition (\ref{eq1:3}) is identically zero. But from (\ref{eq2:3}) it follows that $T_k(t) \equiv 0$ for all $k \geq 1$. Taking into account the
definition of $T_k(t)$ and the completeness of the system $\{v_k\}$, we obtain $u(t) \equiv 0$.

Let us prove that the solution $u(t)$ satisfies the conditions of the Definition \ref{def1}. Define by $S_n(t)$ the partial sum  of the series (\ref{eq2:1}). Then
$$
||A^\sigma S_n(t)||^2=\left|\left|\sum\limits_{k=1}^n \lambda^\sigma[\varphi_k E_{\rho}(-\lambda_k^\sigma t^\rho)]v_k\right|\right|^2.
$$
Apply  Parseval's equality to get
$$
||A^\sigma S_n(t)||^2=\sum_{k=1}^n\left|\lambda_k^\sigma\varphi_kE_{\rho}(-\lambda_k^\sigma t^\rho) \right|^2.
$$
Now using estimate (\ref{eq1:2}), we have
$$
||A^\sigma S_n(t)||^2 \leq C\sum_{k=1}^n|\varphi_k|^2\left|\frac{\lambda_k^\sigma}{1+\lambda_k^\sigma t^\rho}\right|^2\leq Ct^{-2\rho}\sum_{k=1}^n|\varphi_k|^2  \leq C
t^{-2\rho}||\varphi||_H.
$$
From this estimate it follows that $A^\sigma u(t)\in C((0,T]; H)$. In addition, from the equation in problem (\ref{eq1:3}) we have $D_t^\rho S_n(t)= -A^\sigma S_n(t)$ and therefore $D_t^\rho u(t)\in C((0,T]; H)$.

Thus Theorem \ref{t1} is proved.

\section{Inverse problems}\label{s:5}

First, we study the inverse problem of determining the order of the fractional derivative. 

Apply equality (\ref{q1}) to get:
$$
|\varphi_1| E_\rho(-\lambda_1^\sigma t_0^\rho)=d_0.
$$
Obviously, for the first inverse problem to have a solution, the number $\frac{d_0}{|\varphi_1|}$ must be in the range of the Mittag-Leffler function $E_{\rho}(-\lambda_1^\sigma t_0 ^{\rho })$. In Lemma \ref{lem1} it is proven that the Mittag-Leffler function decreases monotonically in $\rho\in[\rho_0,1)$. Therefore, $E_{\rho_0}(- t_0^{\rho_0})>e^{-\lambda_1^\sigma t_0}$. This implies the necessity and sufficiency of inequality (\ref{ty1}). If this condition is met, then, due to the monotonicity of the Mittag-Leffler function, the unknown parameter $\rho$ exists as a solution to the equation (\ref{q1}) and it is unique.

Theorem \ref{t2} is proven.

Now we turn to the study of the inverse problem of simultaneous determination of the orders of the fractional derivative and the degree of the Laplace operator and prove the Theorem \ref{t3}. 

From the additional conditions (\ref{q1}) and (\ref{q2}) we have:
\begin{equation}\label{1}
	\left\{
	\begin{aligned}
		&\varphi_1 E_\rho (-\lambda_1^\sigma t_0^\rho)-d_0=0, \\
		&\varphi_1 E_\rho (-\lambda_1^\sigma t_1^\rho)-d_1=0. \\
	\end{aligned}
	\right.
\end{equation}
This is a system of equations for finding the unknown parameters $\rho$ and $\sigma$.

Let us denote by $\mathbb{D}$ the following determinant:  
\begin{equation}\label{2}
	\mathbb{D}=\left|
	\begin{aligned}
		&\frac{\partial}{\partial\rho}E_\rho (-\lambda_1^\sigma t_0^\rho) \quad \frac{\partial}{\partial\sigma}E_\rho (-\lambda_1^\sigma t_0^\rho) \\
		&\frac{\partial}{\partial\rho}E_\rho (-\lambda_1^\sigma t_1^\rho) \quad \frac{\partial}{\partial\sigma}E_\rho (-\lambda_1^\sigma t_1^\rho) \\
	\end{aligned}
	\right|.
\end{equation}
According to the theorem on the existence of a solution to a system of equations (see \cite{IP}, p. 581), in order for the system (\ref{1}) to have a unique solution, it is necessary and sufficient that $\mathbb{D}$ be different from zero for all $\rho\in [\rho_0, \rho_1]$ and $\sigma>0$.

Let us calculate the necessary derivatives. For this we use equality (\ref{eq3:1}) by putting the values $t_0$ and $t_1$ respectively to the functions $p(\rho,\sigma,t)$ and $q(\rho,\sigma,t)$:
$$
\frac{\partial}{\partial\rho}E_\rho (-\lambda_1^\sigma t_0^\rho)=-\frac{\ln t_0-\Psi(1-\rho)}{\Gamma(1-\rho)\lambda_1^\sigma t_0^\rho}+\frac{\partial}{\partial\rho}q(\rho,\sigma,t_0),
$$

$$
\frac{\partial}{\partial\sigma}E_\rho (-\lambda_1^\sigma t_0^\rho)=-\frac{\ln \lambda_1}{\Gamma(1-\rho)\lambda_1^\sigma t_0^\rho}+\frac{\partial}{\partial\sigma}q(\rho,\sigma,t_0),
$$

$$
\frac{\partial}{\partial\rho}E_\rho (-\lambda_1^\sigma t_1^\rho)=-\frac{\ln t_1-\Psi(1-\rho)}{\Gamma(1-\rho)\lambda_1^\sigma t_1^\rho}+\frac{\partial}{\partial\rho}q(\rho,\sigma,t_1),
$$

$$
\frac{\partial}{\partial\sigma}E_\rho (-\lambda_1^\sigma t_1^\rho)=-\frac{\ln \lambda_1}{\Gamma(1-\rho)\lambda_1^\sigma t_1^\rho}+\frac{\partial}{\partial\sigma}q(\rho,\sigma,t_1).
$$

Now we have
$$
\mathbb{D} = \frac{\ln\lambda_1(\ln t_0-\Psi(1-\rho))}{\Gamma^2(1-\rho)\lambda_1^{2\sigma}t_0^\rho t_1^\rho}-\frac{\ln \lambda_1(\ln t_1-\Psi(1-\rho))}{\Gamma^2(1-\rho)\lambda_1^{2\sigma}t_0^\rho t_1^\rho}-\frac{\ln\lambda_1}{\Gamma(1-\rho)\lambda_1^\sigma t_1^\rho}\frac{\partial}{\partial\rho}q(\rho,\sigma,t_0)-
$$
$$
-\frac{\ln t_0-\Psi(1-\rho)}{\Gamma(1-\rho)\lambda_1^\sigma t_0^\rho}\frac{\partial}{\partial\sigma}q(\rho,\sigma,t_1)+\frac{\ln\lambda_1}{\Gamma(1-\rho)\lambda_1^\sigma t_0^\rho}\frac{\partial}{\partial\rho}q(\rho,\sigma,t_1)+
$$
$$
+\frac{\ln t_1-\Psi(1-\rho)}{\Gamma(1-\rho)\lambda_1^\sigma t_1^\rho}\frac{\partial}{\partial\sigma}q(\rho,\sigma,t_0)+\frac{\partial}{\partial\rho}q(\rho,\sigma,t_0)\frac{\partial}{\partial\sigma}q(\rho,\sigma,t_1)-\frac{\partial}{\partial\sigma}q(\rho,\sigma,t_0)\frac{\partial}{\partial\rho}q(\rho,\sigma,t_1),
$$
or
$$
\mathbb{D} = \frac{\ln \lambda_1(\ln t_0-\ln t_1)}{\Gamma^2(1-\rho)\lambda_1^{2\sigma}t_0^\rho t_1^\rho}-\frac{\ln\lambda_1}{\Gamma(1-\rho)\lambda_1^\sigma t_1^\rho}\frac{\partial}{\partial\rho}q(\rho,\sigma,t_0)-
$$
$$
-\frac{\ln t_0-\Psi(1-\rho)}{\Gamma(1-\rho)\lambda_1^\sigma t_0^\rho}\frac{\partial}{\partial\sigma}q(\rho,\sigma,t_1)+\frac{\ln\lambda_1}{\Gamma(1-\rho)\lambda_1^\sigma t_0^\rho}\frac{\partial}{\partial\rho}q(\rho,\sigma,t_1)+
$$
$$
+\frac{\ln t_1-\Psi(1-\rho)}{\Gamma(1-\rho)\lambda_1^\sigma t_1^\rho}\frac{\partial}{\partial\sigma}q(\rho,\sigma,t_0)+\frac{\partial}{\partial\rho}q(\rho,\sigma,t_0)\frac{\partial}{\partial\sigma}q(\rho,\sigma,t_1)-\frac{\partial}{\partial\sigma}q(\rho,\sigma,t_0)\frac{\partial}{\partial\rho}q(\rho,\sigma,t_1).
$$
Without loss of generality, we assume that $t_0>t_1$.
Consider two cases: the first case $\lambda_1>1$ and the second case $0<\lambda_1<1$. 

If $\lambda_1>1$, then apply the estimates (\ref{es1}), (\ref{es2}) and (\ref{es4}) to get:

$$
\mathbb{D}\geq
\frac{\ln\lambda_1(\ln t_0-\ln t_1)}{\Gamma^2(1-\rho)\lambda_1^{2\sigma}t_0^\rho t_1^\rho}-C_1\frac{\ln\lambda_1(1/\rho+\ln t_0)}{\Gamma(1-\rho)\lambda_1^{3\sigma} t_0^{2\rho} t_1^\rho}-\frac{C_3\ln\lambda_1}{\lambda_1^{3\sigma}t_0^\rho t_1^\rho}-
$$
$$
-C_2\frac{\ln\lambda_1(1/\rho+\ln t_1)}{\Gamma(1-\rho)\lambda_1^{3\sigma}t_0^\rho t_1^{2\rho}}-\frac{C_3\ln\lambda_1}{\lambda_1^{3\sigma}t_0^\rho t_1^\rho} -C_1C_3\frac{\ln\lambda_1(1/\rho+\ln t_0)}{\lambda_1^{4\sigma} t_0^{2\rho}t_1^\rho}-
$$
$$
-C_2C_3\frac{\ln\lambda_1(1/\rho+\ln t_1)}{\lambda_1^{4\sigma} t_0^\rho t_1^{2\rho}},
$$
where $C_1$ and $C_2$ are defined in Lemma \ref{lem1} and $C_3$ is defined in Lemma \ref{lem4}.
From this, we get:
$$
\mathbb{D}\geq
\frac{\ln\lambda_1}{\lambda_1^{2\sigma}t_0^\rho t_1^\rho}\left(\frac{\ln t_0-\ln t_1}{\Gamma^2(1-\rho)}-C_1\frac{1/\rho+\ln t_0}{\Gamma(1-\rho)\lambda_1^\sigma t_0^\rho}-C_2\frac{1/\rho+\ln t_1}{\Gamma(1-\rho)\lambda_1^\sigma t_1^\rho}-\right.
$$
$$
\left. -C_1'\frac{1/\rho+\ln t_0}{\lambda_1^{2\sigma} t_0^\rho}-C_2'\frac{1/\rho+\ln t_1}{\lambda_1^{2\sigma} t_1^\rho}-\frac{2C_3}{\lambda_1^\sigma}\right)=\mathbb{D}_1.
$$
where $C_1'=C_1C_3,\quad C_2'=C_2C_3$.

Using $\Gamma(1-\rho)\geq 1, \quad \lambda_1^{2\sigma}\geq \lambda_1^{\sigma}$ and equality (\ref{eq3:2}) we obtain:

$$
\mathbb{D}_1\geq \frac{\ln\lambda_1}{\lambda_1^{2\sigma}t_0^\rho t_1^\rho}\left(\frac{(1-\rho)^2(\ln t_0-\ln t_1)}{\Gamma^2(2-\rho)}-C_1''\frac{(1/\rho+\ln t_0)}{\lambda_1^\sigma t_0^\rho}-\right.
$$
$$
\left. -C_2''\frac{(1/\rho+\ln t_1)}{\lambda_1^\sigma t_1^\rho}-\frac{2C_3}{\lambda_1^\sigma}\right)=\mathbb{D}_2.
$$
where $C_1''=C_1'+1,\quad C_2''=C_2+1$.

By inequality $\Gamma(2-\rho)\leq 1$ and choosing $t_0$ and $t_1$ sufficiently large, 
in particular
$$
t_0^{\rho_0}>C_1''(1/{\rho_0}+\ln{t_0}),\quad t_1^{\rho_0}>C_2''(1/{\rho_0}+\ln{t_1}),
$$
we get:
$$
\mathbb{D}_2\geq \frac{\ln\lambda_1}{\lambda_1^{2\sigma}t_0^\rho t_1^\rho}\left((1-\rho)^2(\ln t_0-\ln t_1)-C \right).
$$
where $C=2(1+C_3)$.

In order for the resulting expression to be positive, the following inequality must be satisfied for $t_0$ and $t_1$:
$$
t_0>t_1e^{\frac{C}{(1-\rho)^2}}.
$$
or
\begin{equation}\label{baho1}
	t_0>t_1e^{\frac{C}{ (1-\rho_1)^2}}.
\end{equation}

If $\lambda_1\in(0,1)$, then
$$
\mathbb{D} = \frac{\ln \lambda_1(\ln t_0-\ln t_1)}{\Gamma^2(1-\rho)\lambda_1^{2\sigma}t_0^\rho t_1^\rho}-\frac{\ln\lambda_1}{\Gamma(1-\rho)\lambda_1^\sigma t_1^\rho}\frac{\partial}{\partial\rho}q(\rho,\sigma,t_0)-
$$
$$
-\frac{\ln t_0-\Psi(1-\rho)}{\Gamma(1-\rho)\lambda_1^\sigma t_0^\rho}\frac{\partial}{\partial\sigma}q(\rho,\sigma,t_1)+\frac{\ln\lambda_1}{\Gamma(1-\rho)\lambda_1^\sigma t_0^\rho}\frac{\partial}{\partial\rho}q(\rho,\sigma,t_1)+
$$
$$
+\frac{\ln t_1-\Psi(1-\rho)}{\Gamma(1-\rho)\lambda_1^\sigma t_1^\rho}\frac{\partial}{\partial\sigma}q(\rho,\sigma,t_0)+\frac{\partial}{\partial\rho}q(\rho,\sigma,t_0)\frac{\partial}{\partial\sigma}q(\rho,\sigma,t_1)-\frac{\partial}{\partial\sigma}q(\rho,\sigma,t_0)\frac{\partial}{\partial\rho}q(\rho,\sigma,t_1).
$$

Apply the estimates (\ref{es1}), (\ref{es2}) and (\ref{es4}) to get:

$$
\mathbb{D}\leq
\frac{\ln\lambda_1(\ln t_0-\ln t_1)}{\Gamma^2(1-\rho)\lambda_1^{2\sigma}t_0^\rho t_1^\rho}-C_1\frac{\ln\lambda_1(1/\rho+\ln t_0)}{\Gamma(1-\rho)\lambda_1^{3\sigma} t_0^{2\rho} t_1^\rho}+\frac{C_3\ln\lambda_1}{\lambda_1^{3\sigma}t_0^\rho t_1^\rho}-
$$
$$
-\frac{C_3\ln\lambda_1}{\lambda_1^{3\sigma}t_0^\rho t_1^\rho} 
-C_2C_3\frac{\ln\lambda_1(1/\rho+\ln t_1)}{\lambda_1^{4\sigma} t_0^\rho t_1^{2\rho}}=
\frac{\ln\lambda_1(\ln t_0-\ln t_1)}{\Gamma^2(1-\rho)\lambda_1^{2\sigma}t_0^\rho t_1^\rho}-
$$
$$
-C_1\frac{\ln\lambda_1(1/\rho+\ln t_0)}{\Gamma(1-\rho)\lambda_1^{3\sigma} t_0^{2\rho} t_1^\rho}-C_2'\frac{\ln\lambda_1(1/\rho+\ln t_1)}{\lambda_1^{4\sigma} t_0^\rho t_1^{2\rho}},
$$
where $C_1$ and $C_2$ are defined in Lemma \ref{lem1}, $C_3$ is defined in Lemma \ref{lem4} and $C_2'=C_2C_3$.
From this, we get:
$$
\mathbb{D}\leq
\frac{\ln\lambda_1}{\lambda_1^{2\sigma}t_0^\rho t_1^\rho}\left(\frac{\ln t_0-\ln t_1}{\Gamma^2(1-\rho)}-C_1\frac{1/\rho+\ln t_0}{\Gamma(1-\rho)\lambda_1^\sigma t_0^\rho}-\frac{C_2'(1/\rho+\ln{t_1})}{\lambda_1^{2\sigma} t_1^\rho}\right)=\mathbb{D}_1.
$$

Apply inequalities $\Gamma(1-\rho)\geq 1, \quad \lambda_1^{2\sigma}\leq \lambda_1^{\sigma}$ and equality (\ref{eq3:2}) to obtain:

$$
\mathbb{D}_1\leq \frac{\ln\lambda_1}{\lambda_1^{2\sigma}t_0^\rho t_1^\rho}\left(\frac{(1-\rho)^2(\ln t_0-\ln t_1)}{\Gamma^2(2-\rho)}-C_1\frac{(1/\rho+\ln t_0)}{\lambda_1^{2\sigma} t_0^\rho}-C_2'\frac{(1/\rho+\ln t_1)}{\lambda_1^{2\sigma} t_1^\rho}\right).
$$
Again by inequality $\Gamma(2-\rho)\leq 1$ and choosing $t_0$ sufficiently large, 
in particular
$$
t_0^\rho>\frac{C_1(1/\rho+\ln{t_0})}{\lambda_1^{2\sigma}},
\quad t_0^\rho>\frac{C_2'(1/\rho+\ln{t_0})}{\lambda_1^{2\sigma}}
$$
or
$$
t_0^{\rho_0}>\frac{C_1(1/\rho+\ln{t_0})}{\lambda_1^{2}}, \quad
t_1^{\rho_0}>\frac{C_2'(1/\rho+\ln{t_0})}{\lambda_1^{2}},
$$
we get:

$$
\mathbb{D}_1\leq \frac{\ln\lambda_1}{\lambda_1^{2\sigma}t_0^\rho t_1^\rho}\left((1-\rho)^2(\ln t_0-\ln t_1)-2\right).
$$

In order for the resulting expression to be negative, the following inequality must be satisfied for $t_0$ and $t_1$:
\begin{equation}\label{baho2}
	t_0>t_1e^{\frac{2}{ (1-\rho_1)^2}}.
\end{equation}

Hence in both cases, one has $\mathbb{D}\neq 0$ for all $\rho\in[\rho_0,\rho_1]$ and $\sigma>0$.

The necessity and sufficiency of
inequality  (\ref{ty2}) can be proved in the same way as for estimate (\ref{ty1}).

Thus, Theorem \ref{t3} is proved.

\section*{Conclusions} When modeling various processes (physical, biological, etc.), specialists increasingly use fractional order equations (see, for example, \cite{Ma}). But when considering such equations, the order of fractional derivatives or the degree of the operator included in the equation is not always known in advance, and, very importantly, there are no tools to measure them. Therefore, over the past 10-15 years, experts have been actively studying inverse problems to determine these unknown parameters..

As noted above, only the next two papers \cite{TU} and \cite{Yamamoto} studied inverse problems for simultaneous determination of two parameters $\rho$ and $\sigma$ in equation (\ref{laplas}). Using additional conditions found in these works, the authors proved the uniqueness of the unknown parameters. It should be emphasized that, using these additional conditions it is almost impossible to prove the existence of a solution to the inverse problem.  We also note that in these papers to prove uniqueness, fairly strict conditions were imposed on the initial function (the initial function must belong to the Sobolev class $L_2^a(\Omega)$, $a>N/2$, and be positive).

In the present paper, new easily verifiable additional conditions are proposed for the inverse problem of determining two parameters $\rho$ and $\sigma$ in equation (\ref{laplas}). As the main results of the work show, these conditions ensure both the existence and uniqueness of the unknown parameters $\rho$ and $\sigma$. It is important to note that when proving these results, it is sufficient that the initial function is from class $L_2(\Omega)$.

\section{Acknowledgments}
The authors are grateful to Sh. Alimov for discussions of these results. The authors acknowledges financial support from the Ministry of Innovative Development of the Republic of Uzbekistan, Grant No F-FA-2021-424.

\bibliographystyle{amsplain}

\end{document}